\documentclass[english,12pt]{amsart}
\usepackage[english]{babel}
\usepackage{amsmath,amssymb,enumerate,amsthm}
\usepackage[latin1]{inputenc}
\usepackage[T1]{fontenc}
\usepackage{cite}
\usepackage{latexsym}
\usepackage{setspace}
\usepackage{bm} 
\usepackage{comment}
\usepackage{array,graphicx,caption}

\newtheorem{lemma}{Lemma}[section]
\newtheorem{theorem}[lemma]{Theorem}

\newcommand{\Bad}{\mathrm{Bad}}
\newcommand{\Tad}{\mathrm{Bad}_{\pmb{\theta}}}
\newcommand{\bi}{\textbf{k}}

\def\={\;=\;}
\def\>{\;>\;}
\def\<{\;<\;}
\def\:{\,:\,}
\def\.={\;\dot{=}\;}

\newcommand{\R}{\mathbb{R}}
\newcommand{\Q}{\mathbb{Q}}
\newcommand{\Z}{\mathbb{Z}}
\newcommand{\N}{\mathbb{N}}

\newcommand{\PP}{\mathcal{P}}

\begin{document}

\begin{Large}
\title{Sets of inhomogeneous linear forms can be not isotropically winning}
\vskip+1cm
\author[N. DYAKOVA]{Natalia Dyakova$^*$}
\thanks{$^*$ Supported by RFBR and RNF grant.}
\address{Department of Mathematics and Mechanics, Moscow State University, Leninskie Gory 1, GZ MGU, 119991 Moscow, Russia}
\email{natalia.stepanova.msu@gmail.com}
\end{Large}

\begin{abstract}
We give an example of irrational  vector $\pmb{\theta} \in \mathbb{R}^2$ such that the set $\Bad_{\pmb{\theta}}  :=  \{(\eta_1,\eta_2):  \inf_{x\in\N} x^{\frac{1}{2}} \max_{i=1,2} \|x \theta_i-\eta_i\|>0\}$  is  not absolutely winning with respect to McMullen's game.
\end{abstract}
\vskip+1.5cm
\begin{abstract}
We give an example of irrational  vector $\pmb{\theta} \in \mathbb{R}^2$ such that the set $\Bad_{\pmb{\theta}}  :=  \{(\eta_1,\eta_2):  \inf_{x\in\N} x^{\frac{1}{2}} \max_{i=1,2} \|x \theta_i-\eta_i\|>0\}$  is  not absolutely winning with respect to McMullen's game.
\end{abstract}
\maketitle
\section{Introduction}
In this paper we consider a problem related to  inhomogeneous Diophantine approximation. Given $\pmb{\theta} =(\theta_1,\theta_2)\in \mathbb{R}^2 $ we study    the set of couples $(\eta_1,\eta_2)\in\mathbb{R}^2$ such that the system of two linear forms
$$
|| x\theta_1 - \eta_1|| , ||x\theta_2 - \eta_2||,
$$
where $||\cdot ||$ stands for the distance to the nearest integer
is badly approximable. We prove a statement complementary to our recent result  from
\cite{BMS}.
We construct $\pmb\theta$ such that the set
$$\Tad :=  \{(\eta_1,\eta_2):  \inf_{x\in\N} x^{\frac{1}{2}} \max_{i=1,2} \|x \theta_i-\eta_i\|>0\}
$$ 
is not isotropically winning.

Our paper is organized as follows. In Section 2 we discuss different games appearing in Diophantine problems. In Section 3
we give a brief survey on inhomogeneous badly approximable systems of linear forms and formulate our main result Theorem 3.1.
Sections 4 and 5 are devoted to some auxiliary observations. In Sections 6, 7, 8 we give a proof for Theorem 3.1.

\section{Schmidt's game and its generalizations}
The following game was introduced by Schmidt in \cite{Sch1}, \cite{Sch2}, \cite{Sch3}. Let $0<\alpha, \beta <1$.
Suppose that two players A and B choose by the turn a nested sequence of closed balls 
$$
B_1 \supset A_1 \supset B_2 \supset A_2 \supset \ldots
$$
with the property that the diameters $|A_i| , |B_i|$ of the balls $A_i, B_i$ satisfy
$$
|A_i| = \alpha |B_i|, \qquad |B_{i+1}| = \beta |A_i| \qquad \forall i = 1,2,3,\ldots
$$
for fixed $0<\alpha, \beta <1$. A set $E\subset \R^n$ is called an $(\alpha,\beta)-winning$ if player A has a strategy which guarantees that intersection $\cap A_i$ meets $E$ regardless of the way B chooses to play.  A set $E\supset \R^n$ is called an $\alpha-winning$ set if it is $(\alpha,\beta)$-winning for all $0<\beta <1$.

There are different modifications of Schmidt's game: strong game and  absolute game introduced in \cite{Mcm}, hyperplane absolute game introduced in \cite{WK}, potential game considered  in \cite{FSU} and some others. 
In our recent paper 
 \cite{BMS} we introduced isotropically winning sets. Let us  describe  here some of these generalizations in more details.

The definition of an absolutely winning set was given by  McMullen  \cite{Mcm}. Consider the following game.
Suppose $A$ and $B$ choose by the turn a sequence of balls $A_i$ and $B_i$ such that the sets
$$
B_{1} \supset (B_1 \backslash A_1) \supset B_2 \supset (B_2 \backslash A_2) \supset B_3 \supset \ldots 
$$
are nested.
For fixed $0<\beta <1/3$ we suppose
$$
|B_{i+1}| \ge \beta |B_i| , \qquad |A_i|\le \beta |B_i|.
$$
We say $E$ is an $absolute$ $winning$ set if for all $\beta \in (0, 1/3),$ player A
has a strategy which guarantees that $\cap B_i$ meets $E$ regardless of how B chooses to play. 
Mcmullen proved that an absolute winning set is $\alpha$-winning for all $\alpha <1/2$.
Several examples of absolute winning sets were exhibited by McMullen in \cite{Mcm}. 
In particular  set of badly approximable numbers in $\R$ is absolutely winning. However the set of simultaneously badly approximable vectors in  $\R^n$ for $ n > 1$
is not absolutely winning.  

In \cite{BMS} another strong variant of winning property was given.
We say that a set $E \subset \mathbb{R}^n$ is {\it isotropically winning} if for each $d \le n$ and for each $d$-dimensional affine subspace $\mathcal{A} \subset \mathbb{R}^n$ the intersection $E \cap \mathcal{A}$ is $1/2$-winning
for Schmidt's game considered as a game in $\mathcal{A} $.
It is clear that an absolute winning set is isotropically winning for each $\alpha \le 1/2$.

\section{Inhomogeneous approximations}
The first important result on inhomogeneous appxohimations in  one-dimensional case is due to
Khintchine.
In \cite{KH1}
he proved that there exists an absolute constant $\gamma$ such that for every $\theta\in \mathbb{R}$ there exists $\eta\in \mathbb{R}$
such that
$$
\inf _{q\in \mathbb{Z} } q ||q\theta - \eta|| > \gamma
.$$
Later  (see \cite{KH2}, \cite{KH3}) he
 proved that for given positive numbers $n, m \in \Z$ there exists a positive constant $\gamma_{nm}$ such that for any $m \times n$ real matrix $\pmb{\theta}$ there exists a vector $\pmb{\eta} \in \R^n$ such that 
$$
\inf_{\pmb{x}\in\Z^m \backslash \{0\}} (|| \pmb{\theta} \pmb{x} - \pmb{\eta}||_{\Z^n})^n ||\pmb{x}||^m >\gamma_{nm}
$$
(here $||\cdot||_{\mathbb{Z}^n}$
stands for the distance to the nearest integral point in $\sup$-norm).
 These results are presented in a wonderful book \cite{Cas} by  Cassels.

Jarnik in \cite{J1}, \cite{J2} proved a generalization of this statement.
Suppose $\psi (t)$ is a function decreasing to zero as $t \rightarrow + \infty$.
 Let $\rho (t)$ be the function inverse to the function $t \mapsto 1/\psi(t)$. Suppose that for all $t > 1$ one has $\psi_{\theta}(t) \le \psi(t)$.
Then there exists a vector $\pmb{\eta} \in \R^n$ such that 
$$
\inf_{\pmb{x}\in\Z^m \backslash \{0\}} (|| \pmb{\theta} \pmb{x} - \pmb{\eta}||_{\Z^n})\cdot \rho(8m\cdot ||\pmb{x}||) >\gamma
$$
with appropriate $\gamma = \gamma(n,m)$.

Denote by
$$
\Bad_{\theta} = \left\{\alpha \in [0,1):  \inf_{q\in\N}q\cdot ||q \theta -\alpha||>0\right\}.
$$
It happened that
the winning property of this inhomogeneous Diophantine set was considered quite recently.
In \cite{T} Tseng showed that $\Bad_{\theta}$ is winning for all real numbers $\theta$ in classical Schmidt's sense.
For the corresponding multidimensional sets
$$
\Bad(n,m)=\left\{\pmb{\theta}\in\mathrm{Mat}_{n\times m}(\R):\, \inf_{q\in\Z^m_{\neq 0}}\max_{1\leq i\leq n}(|q|^{\frac{m}{n}}\|\pmb{\theta}_i(q)\|)>0\right\}.
$$
the winning 
property  shown for example in 
  \cite{ET} and \cite{NGM}.

  Further generalizations deal with the {\it twisted} sets
$$
\Bad(i,j)=\left\{(\theta_1,\theta_2)\in\R^2:\, \inf_{q\in\N}\max(q^i\|q\theta_1\|, q^j\|q\theta_2\|)>0\right\},
$$
where $i,j$ are real positive numbers satisfying $i+j=1$,
introduced by Schmidt. 
 An in \cite{An} proved that $\Bad(i,j)$ is winning for the standard Schmidt game.
 In higher dimension, we fix an $n$-tuple $\bi=(k_1,\ldots,k_n)$ of real numbers satisfying
\begin{equation}\label{weight}
k_1,\ldots,k_n>0\qquad\mbox{and}\qquad \sum_{i=1}^n k_i=1,
\end{equation}
and define
$$
\Bad(\bi,n,m)=\left\{\pmb{\theta}\in\mathrm{Mat}_{n\times m}(\R):\, \inf_{q\in\Z^m_{\neq 0}}\max_{1\leq i\leq n}(|q|^{mk_i}\|\pmb{\theta}_i(q)\|)>0\right\}.
$$
Here, $|\cdot|$ denotes the supremum norm, $\pmb{\theta}=(\pmb{\theta}_{ij})$ and $\pmb{\theta}_i(q)$ is the product of the $i$-th line of $\pmb{\theta}$ with the vector $q$, i.e.
$$
\pmb{\theta}_i(q)=\sum^m_{j=1} q_j \pmb{\theta}_{ij}.
$$ 
In twisted setting much less is known. 
In particular up to now winning property of the set $
\Bad(\bi,n,m)$ in dimension greater that two is not proved.

Given $\pmb{\theta}\in\mathrm{Mat}_{n\times m}(\R)$, we define
$$
\Tad(\bi,n,m)=\left\{x\in\R^n:\, \inf_{q\in\Z^m_{\neq 0}}\max_{1\leq i\leq n}(|q|^{mk_i}\|\pmb{\theta}_i(q)-x_i\|)>0\right\}.
$$
Harrap and Moshchevitin in \cite{HM} showed that this set is winning provided that $\pmb{\theta}\in\Bad(\bi,n,m)$.
In \cite{BMS} 
it was proved that
if we suppose that
 $\pmb{\theta}\in\Bad(\bi,n,m)$
 the set $
\Tad(\bi,n,m)$ is isotropically winning\footnote{In fact the approach from \cite{BMS} 
gives a little bit more.
Instead of property that for any subspace $\mathcal{A}$  the intersection $E \cap \mathcal{A}$  is 1/2 winning  in $\mathcal{A}$ one can see that it is $\alpha$-winning for all $\alpha \in (0,1/2]$. It is not completely clear for the author if these two properties are equivalent. (For a closely related problem see \cite{Dremov}.)}.

We should note that even in the case
  $n=2, m=1$
  it is not known if the set
  $\Tad(\bi,2,1)$ is $\alpha$-winning for some positive $\alpha$ without the condition
   $\pmb{\theta}\in\Bad(\bi,2,1)$.
  .

In this article we show that the condition for $\pmb{\theta}$ be from $\Bad(\bi,n,m)$ is essential for the isotropically winning property, and prove the following theorem.
\begin{theorem}\label{th}

There exists a vector $\pmb{\theta}= (\theta_1,\theta_2)$ such that

1) $ 1,\theta_1,\theta_2$ are lineary independent over $\Z$;

2) $\Tad :=  \{(\eta_1,\eta_2):  \inf_{x\in\N} x^{\frac{1}{2}} \max_{i=1,2} \|x \theta_i-\eta_i\|>0\}$ is not isotropically winning.
\end{theorem}

\section{Some more remarks}
In the sequel 
$
\pmb{x} = (x_0, x_1,x_2)
$ is a vector in $\mathbb{R}^3$,
$| \cdot |$  stands  for the Euclidean norm of the vector, and by $(\bm{w}, \bm{t})$ we denote the
innerproduct of vectors $\bm{w} $ and $\bm{t}.$

The proof of the Theorem \ref{th} we will give in Section \ref{indep}.
There we will construct a special $\pmb{\theta}$ and one-dimensional affine subspace $\PP $ such that  $\pmb{\theta} \in \PP$ and for the segment $\mathcal{D} = \PP \cap \{|\pmb{z} -\pmb{\theta}|\le 1\}$  one has $\mathcal{D} \cap \Tad = \emptyset.$
Moreover given arbitrary positive function $\omega(t)$ monotonically (slowly) increasing to infinity we can ensure that for all $\pmb{\eta} = (\eta_1,\eta_2) \in \mathcal{D}$ there exist infinitely many $x \in \Z$ such that 
$$
\max_{i=1,2} \|x \theta_i-\eta_i\|< \dfrac{\omega(x)}{x}.
$$
To explain the construction of the proof it is useful to consider the case when $\theta_1, \theta_2,1$ are linear dependent. This case we will discuss in the Section \ref{dep}.

{\bf Remark 1.}
 From the result of the paper \cite{BMS}  it follows that the vector $\pmb{\theta}$ constructed in Theorem \ref{th}
 does not belong to the set 
 $$\Bad = \{(\theta_1,\theta_2) | \inf_{x \in \N} x^{1/2} \max (||\theta_1x||,||\theta_2x||)>0\}.$$

{\bf Remark 2.}
Let $\pmb{\theta} = \left(\dfrac{a_{1}}{q}, \dfrac{a_{2}}{q}\right)$ be rational.
Let $\pmb{\eta} = (\eta_1,\eta_2) \notin \dfrac{1}{q} \cdot \Z^2$, then for any $x \in \Z$,
$$
\max_{i=1,2} \left|\left|x \dfrac{a_{i}}{q}-\eta_i\right|\right| \ge {\rm dist} \left(\pmb{\eta},\dfrac{1}{q} \cdot \Z^2\right)> 0 .
$$
So the set $$\mathcal{B} = \left\{ \pmb{\eta} :  \inf_{x \in \Z} \max_{i=1,2} \left|\left|x \dfrac{a_{i}}{q}-\eta_i\right|\right| > 0 \right\} $$
contains $\R^2 \setminus \dfrac{1}{q} \cdot \Z^2 $ and is trivially winning.
It is clear, that for any one-dimensional affine subspace $\ell$ we have $\mathcal{B} \cap \ell \supset \left(\R^2 \setminus \dfrac{1}{q} \cdot \Z^2 \right)\cap \ell$.
So obviously $\mathcal{B} \cap \ell$ is also winning in $\ell.$
\section{linear dependent case}\label{dep}
Let $1,\theta_1,\theta_2$ be linear dependent and at least one of $\theta_j$ is irrational.  It means that there exists $\bm{z} = (z_0,z_1,z_2) \in \Z^3$ such that $(\bm{z} , \bm{\theta}) = 0.$
Let us consider two-dimensional rational subspace $$\pi = \{\pmb{x} \in R^3 : (\pmb{x},\pmb{z}) =0\},$$ so $\bm{\theta} \in \pi$.

Let us define one-dimensional subspace $\PP = \{ (x_1, x_2) : (1,x_1,x_2) \in \pi\} \subset \R^2$ .

We will prove that there exists a constant $\gamma$ such that for any 
$
\eta = (\eta_1,\eta_2) \in \PP 
$
the inequality 
$$
\qquad \max_{i=1,2} ||\theta_i x - \eta_i|| < \frac{\gamma}{x}
$$
has  infinitely many solutions in $x \in \N$.
(This statement is similar to Chebushev's theorem \cite[Theorem 24, Chapter 2]{CH})

Denote by $\Lambda = \pi \cap \Z^3$ the integer lattice with the determinant $ d : = \det \Lambda = |\pmb{z}|.$ 
Denote by $\{\pmb{g}_{\nu} = (q_{\nu}, a_{1 \nu}, a_{2 \nu})\}_{\nu=1,2,3 \ldots} \subset \Lambda$ the sequence of the best approximations of $\pmb{\theta}$ by the lattice $\Lambda$ and corresponding parallelograms

$$
\Pi_{\nu} = \{ \pmb{x} = (x_0,x_1,x_2) \in \pi : 0\le x_0 \le q_{\nu} : $$
$${\rm dist} (\pmb{x}, l(\pmb{\theta})) \le {\rm dist} (\pmb{g}_{\nu -1}, l(\pmb{\theta})) \}
$$
which contains a fundamental domain of the two-dimensional $\Lambda.$
Obviously, vol $\Pi_{\nu} \le 4d$.
So 
\begin{equation}\label{1}
 {\rm dist} (\pmb{g}_{\nu -1}, l(\theta)) \ll \dfrac{d}{q_{\nu}}
\end{equation}
with an absolute constant in the sign $\ll$.
 It is clear that for any point $\bm{\eta} \in \pi$ the shift  $  \bm{\eta} + \Pi_{\nu}$ contains a point of $\Lambda$.

For any $\bm{\eta} = (\eta_1,\eta_2) \in \PP$ and for any positive integer $\nu$ the planar domain $\overline{\pmb{\eta}}+ \Pi_{\nu}$, $\overline{\pmb{\eta}} = (1, -\eta_1, -\eta_2)$ contains an integer point
$\bm{y} = (x,y_1,y_2) \in \Lambda$.

It is clear that 
\begin{equation}\label{2}
1\le x \le 1 + q_{\nu}
\end{equation}
and 
$$
\max_{i=1,2} 
 || \theta_i x- \eta_i||  \ll {\rm dist} (\pmb{y},l(\pmb{\theta}) + \overline{\pmb{\eta}}) \ll  {\rm dist} (l(\pmb{\theta}),\pmb{g}_{\nu -1}),
$$
and by \eqref{1} 
\begin{equation}\label{3}
\max_{i=1,2} || \theta_i x- \eta_i|| \ll \dfrac{d}{q_{\nu}}
\end{equation}

From \eqref{2}, \eqref{3} it follows that the inequality
$$ \max_{i=1,2} || \theta_i x- \eta_i|| \ll \dfrac{d}{x}.$$
has infinitely many solutions
and everything is proved.

\section{Inductive construction of integer points}\label{indep}
Let $\omega(t)$ be arbitrary positive function monotonically (slowly) increasing to infinity. Here we describe the inductive construction of integer points
$\pmb{z}_{\nu}=(q_{\nu},z_{1 \nu },z_{2 \nu })$.
The base of the induction process is trivial. One can take an arbitrary primitive pair of integer vectors that can be completed to a basis of $\mathbb{Z}^3$.

Suppose that we have two primitive integer vectors 
$$
\pmb{z}_{\nu-1}=(q_{\nu-1},z_{1 \, \nu-1 },z_{2 \, \nu-1 })\in \Z^3
$$
$$
\pmb{z}_{\nu}=(q_{\nu},z_{1 \, \nu},z_{2 \, \nu})\in \Z^3.
$$
Now we explain how to construct the next integer vector $\pmb{z}_{\nu+1}.$

We consider two-dimensional subspace
$$\pi_{\nu} = <\pmb{z}_{\nu-1},\pmb{z}_{\nu}>_{\R}. $$
As the pair of vectors $\pmb{z}_{\nu-1}$ and $\pmb{z}_{\nu}$  is primitive, the lattice 
$$\Lambda_{\nu} := <\pmb{z}_{\nu-1},\pmb{z}_{\nu}>_{\Z}=  \pi_{\nu}\cap \Z^3 .$$
By $d_{\nu} = \det \Lambda_{\nu}$ we denote the two-dimensional fundamental volume of the lattice $\Lambda_{\nu}.$
Now we define vector $\pmb{n_{\nu}} = (n_{0 \nu},n_{1\nu}, n_{2\nu})\in \R^3$ from the conditions
$$\pi_{\nu}= \{\pmb{x} \in \R^3 : (\pmb{x},\pmb{n_{\nu}}) =0 \}, \quad |\pmb{n}_{\nu}|=1.$$
Put 
\begin{equation}\label{sigma}
 \sigma_{\nu}= {\rm dist} (\pmb{z}_{\nu -1}, l(\pmb{z}_{\nu}))
\end{equation}
Obviously, $|\pmb{z}_{\nu}| \asymp q_{\nu}$ and
\begin{equation}\label{asymp}
\sigma_{\nu}\asymp \dfrac{d_{\nu}}{q_{\nu}}.
\end{equation}
We define a vector $\pmb{e}_{\nu}$
from the conditions
\begin{equation}\label{e}
\pmb{e}_{\nu} \in \pi_{\nu} , \qquad |\pmb{e}_{\nu}| = 1, \qquad (\pmb{e}_{\nu}, \pmb{z}_{\nu}) = 0
\end{equation}
so $\pmb{e}_{\nu}$ is parallel to $\pi_{\nu}$ and orthogonal to $\pmb{z}_{\nu}$.

Define the rectangle $$\Pi_{\nu}
= \{ \pmb{x} = (x_0,x_1,x_2)  : \pmb{x} = t\pmb{z}_{\nu} + r\pmb{e}_{\nu}, \,\, 0 \le t\le |\pmb{z}_{\nu}|, \,\, |r| \le \sigma_{\nu}  \}.$$
It is clear that rectangle $\Pi_{\nu} \subset \pi_{\nu}$ contains a fundamental domain of the lattice $\Lambda_{\nu}.$
We need two axillary vectors $\pmb{z}^a_{\nu}$ and $\pmb{z}^b_{\nu}$ defined as
$$\pmb{z}^a_{\nu} = \pmb{z}_{\nu} + a_{\nu} \pmb{e}_{\nu}$$ 
$$
\pmb{z}^{b}_{\nu} = \pmb{z}^a_{\nu} + b_{\nu} \pmb{n}_{\nu}
$$
where positive $a_{\nu}$ is chosen in such a way 
\begin{equation}\label{a}
a_{\nu} d_{\nu}^2 \le \nu^{-1} \omega \left( \dfrac{q_{\nu}^2}{d_{\nu}^2} \cdot \dfrac{1}{a_{\nu}} \right) 
\end{equation}
and
\begin{equation}\label{b}
b_{\nu} = a_{\nu} {\rm min} \left(1, \dfrac{d_{\nu}}{q_{\nu}}\right).
\end{equation}

From the construction it follows that 
\begin{equation}\label{|Z|}
|\pmb{z}_{\nu}^a| \asymp |\pmb{z}_{\nu}^b| \asymp |\pmb{z}_{\nu}| \asymp q_{\nu}.  
\end{equation}
Integer lattice $\Z^3$ splits into levels with respect to the two-dimensional sublattice $\Lambda_{\nu}$ in such a way that
$$\Z^3 = \bigsqcup_{i\in \Z} \Lambda_{\nu,i},$$
where $\Lambda_{\nu,j}=\Lambda_{\nu} + j \pmb{z'},  j\in \Z$ and integer vector $\pmb{z'}$ completes the couple $\pmb{z}_{\nu-1}, \pmb{z}_{\nu}$ to the basis in $\Z^3$.
We consider affine subspace 
$\pi^1_{\nu}= \pi_{\nu} + \pmb{z}' \supset \Lambda_{\nu,1},$ which is parallel to $\pi_{\nu}$.
It is clear that ${\rm dist} (\pi_{\nu},\pi^1_{\nu}) = \frac{1}{d_{\nu}}$.

We need to determine the next integer point $z_{\nu+1}$.
Denote by $\mathfrak{P}$ central projection with center 0 onto the affine subspace
 $\pi^1_{\nu}$.
  We consider triangle $\Delta$ with the vertices $z_{\nu},z^a_{\nu},z^{b}_{\nu}$ and its image $\mathfrak{P} \Delta$ under the projection $\mathfrak{P}$.
   Define \begin{equation}\label{Z}
  \pmb{Z} = \mathfrak{P} {\pmb{z}^b_{\nu}}.
   \end{equation}
    One can see that 
  \begin{equation}\label{|zz|}
  |\pmb{Z}| \asymp \dfrac{q_{\nu}}{d_{\nu} b_{\nu}}.
  \end{equation}
  \includegraphics[width=10cm,height=5cm]{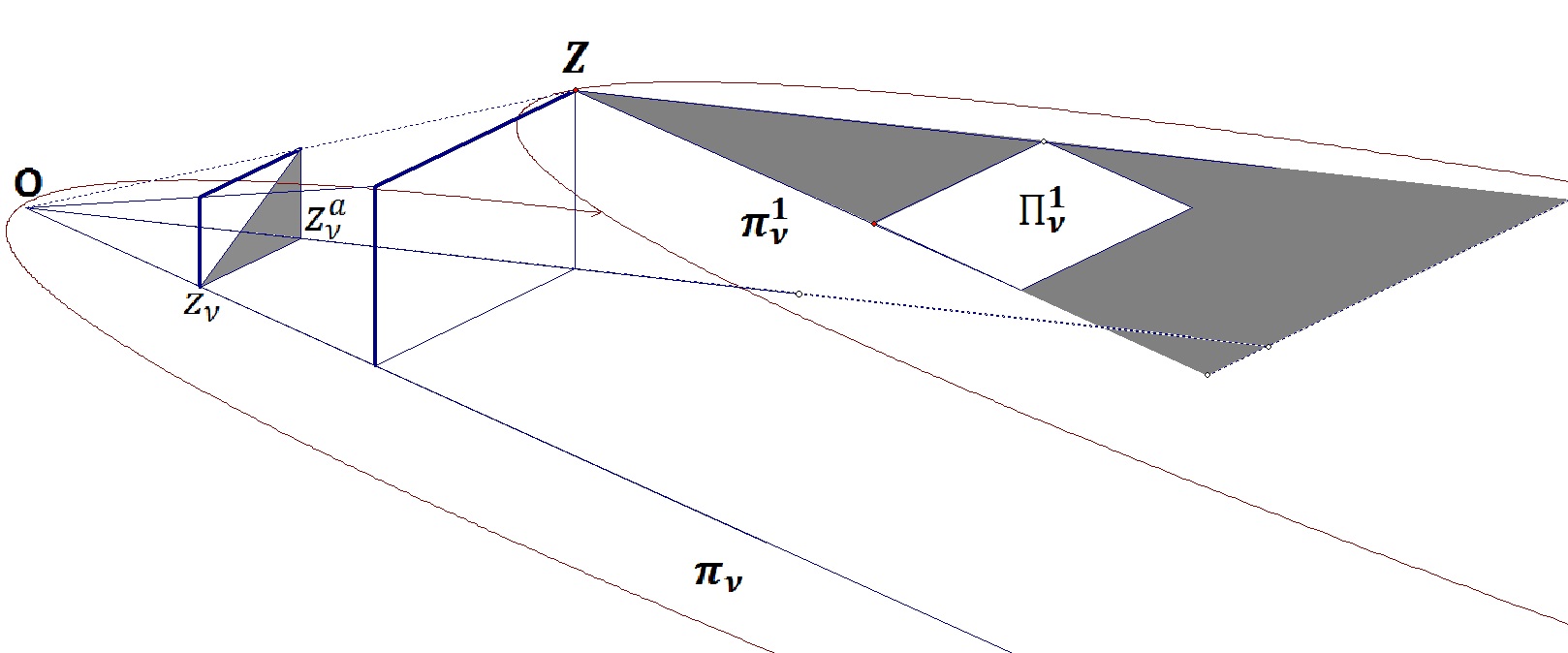} 
  \captionof{figure}{The central projection $ \mathfrak{P}$.}\label{fig:intab}
  
   
Define rays   
$$ \mathcal{R}_1 = \{\pmb{z} = \pmb{Z} + t \pmb{z}_{\nu}, \quad t\ge 0\} \quad \text{and} \quad \mathcal{R}_2 = \{\pmb{z} = \pmb{Z} + t \pmb{z}^a_{\nu}, \quad t\ge 0\}.$$
It is clear that $\mathcal{R}_1 \cap \mathcal{R}_2 = \{\pmb{Z}\}$ and $\mathcal{R}_1, \mathcal{R}_2 \subset \pi^1_{\nu}$.
Moreover, the whole convex angle bounded by rays $\mathcal{R}_1, \mathcal{R}_2$ form the image of the triangle $\Delta$ under the projection $\mathfrak{P}$:
$$
\mathfrak{P}\Delta = {\rm conv}\, (\mathcal{R}_1 \cup \mathcal{R}_2) .
$$
Affine subspace $\pi_{\nu}^1$ contains the affine lattice $\Lambda_{\nu}^1 = \Lambda_{\nu} + \pmb{z}'$ which is congruent to the lattice $\Lambda_{\nu}.$
Thus, for any $\pmb{\zeta} \in \pi_{\nu}^1$ the shift $\Pi_{\nu}+\pmb{\zeta}$ contains an integer point from $\Lambda_{\nu}^1.$

Put 
\begin{equation}\label{tau}
\tau_{\nu} = \dfrac{2\sigma_{\nu} |\pmb{z}_{\nu}|}{a_{\nu}},
\end{equation}
consider the point 
$$\pmb{\zeta}_{\nu} =  \pmb{Z}+ \tau_{\nu} \pmb{z}_{\nu} +\sigma_{\nu}\pmb{e}_{\nu}\in \pi_{\nu}^1. $$
and the rectangle
$$
\Pi_{\nu}^1 = \Pi_{\nu} +\pmb{\zeta}_{\nu} \subset \pi_{\nu}^1.
$$
It is clear that 
$$
\Pi_{\nu}^1 \subset \mathfrak{P} \Delta$$
(here $\pmb{Z}$ was defined in \eqref{Z}, $\pmb{e}_{\nu}$ was defined in \eqref{e} and parameters $\sigma_{\nu}, \tau_{\nu}$ come from \eqref{sigma} and \eqref{tau}).

Now we take the integer point
$$\bm{z}_{\nu+1} = (q_{\nu+1},z_{1 \, \nu+1 },z_{2 \, \nu+1 }) \in \Lambda_{\nu}^1 \cap \Pi_{\nu}^1.$$
From the construction it follows that 
$$q_{\nu+1} \asymp |\pmb{z}_{\nu+1}| \asymp  |\pmb{Z}|+ \tau_{\nu}  |\pmb{z}_{\nu}| + |\pmb{z}_{\nu}| 
$$
$$
\asymp q_{\nu} \left( 1+ \dfrac{1}{d_{\nu} b_{\nu}} + \dfrac{\sigma_{\nu}}{a_{\nu}} \right)  \asymp
 q_{\nu} \left( 1+ \dfrac{1}{d_{\nu} b_{\nu}} \right)+ \dfrac{d_{\nu}}{a_{\nu}} \asymp \dfrac{q_{\nu}}{d_{\nu} b_{\nu}}.  
$$(Here we use (\ref{asymp},\ref{b}, \ref{|Z|}, \ref{|zz|},  \ref{tau}).)
From \eqref{b} we see that 
\begin{equation}\label{q+1}
q_{\nu+1} \gg \left( \dfrac{q_{\nu}}{d_{\nu}}\right)^2 \dfrac{1}{a_{\nu}}.
\end{equation}
Now we are able to define the next  two-dimensional lattice
$$
\Lambda_{\nu+1} = <\pmb{z}_{\nu},\pmb{z}_{\nu+1}>_{\Z}.
$$
Let $d_{\nu+1}$ be its fundamental volume. We will estimate the value of $d_{\nu+1}$ taking into account \eqref{b}
\begin{equation}\label{d}
d_{\nu+1} \ll q_{\nu} \cdot {\rm dist} (\pmb{z}_{\nu+1}, l(\pmb{z}_{\nu})) \ll \dfrac{q_{\nu}}{d_{\nu}} \cdot \dfrac{a_{\nu}}{b_{\nu}} \ll \left(\dfrac{q_{\nu}}{d_{\nu}} \right)^2 \ll q_{\nu}^2.
\end{equation}
From \eqref{q+1} and \eqref{d} we deduce that
$$
d_{\nu+1}  \ll a_{\nu}  d_{\nu}^2 q_{\nu+1}.
$$
By the choice of $a_{\nu}$ (by formula \eqref{a}) we have
\begin{equation}\label{7}
d_{\nu+1} \le \dfrac{\omega(q_{\nu+1})}{\nu}.
\end{equation}
\section{The vector $\pmb{\theta}$}

Now we define
$$
\pmb{\theta}_{\nu} = (\theta_{1 \nu},\theta_{2 \nu}), \quad \theta_{j \nu} = \dfrac{q_{j \nu}}{q_{\nu}}.
$$
We consider the angles between the successive vectors $\pmb{n}_{\nu}$ and $\pmb{n}_{\nu+1}$:
$$
\alpha_{\nu} ={\rm angle} (\pmb{n}_{\nu}, \pmb{n}_{\nu+1}) \asymp \tan {\rm angle} (\pmb{n}_{\nu}, \pmb{n}_{\nu+1}) 
$$ 
 \includegraphics[width=15cm,height=5cm]{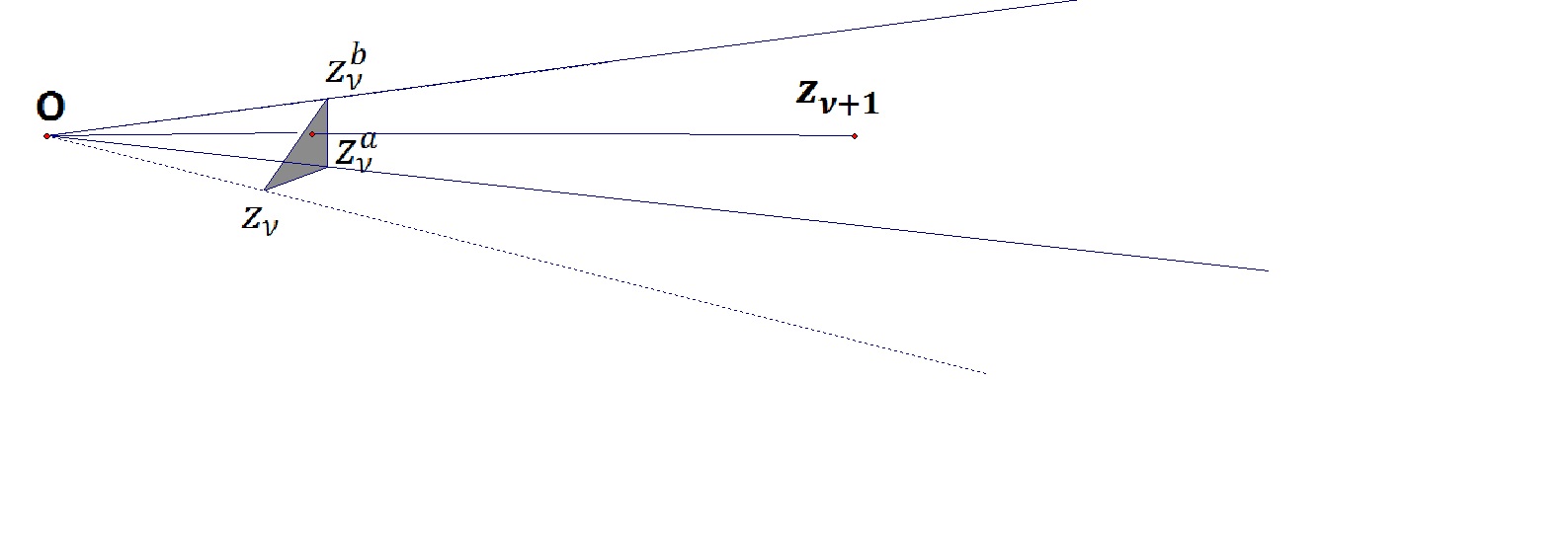} 
 \captionof{figure}{The vector $\pmb{z}_{\nu+1}$ intersects the interior of the triangle    
$\Delta = \pmb{z}_{\nu} \pmb{z}_{\nu}^a\pmb{z}_{\nu}^b$.}\label{fig:intab}


As $\pmb{z}_{\nu+1} \in \mathfrak{P} \Delta$ (see Fig. 2), we have 
$$
\tan  {\rm angle} (\pmb{n}_{\nu}, \pmb{n}_{\nu+1}) \le \dfrac{b_{\nu}}{a_{\nu}},
$$
and so
\begin{equation}\label{alpha}
\alpha_{\nu} \ll \dfrac{b_{\nu}}{a_{\nu}}.
\end{equation}
As $\pmb{z}_{\nu+1} \in \mathfrak{P} \Delta$ we have 
\begin{equation}\label{theta}
|\pmb{\theta}_{\nu} - \pmb{\theta}_{\nu+1}| \ll \dfrac{\sqrt{a_{\nu}^2 + b_{\nu}^2}}{q_{\nu}} \ll \dfrac{a_{\nu}}{q_{\nu}},
\end{equation}
by the same argument.
There exist limits $$\lim_{\nu \rightarrow \infty} \pmb{\theta}_{\nu} = \pmb{\theta} = (\theta_1,\theta_2) \quad {\rm and} \quad \lim_{\nu \rightarrow \infty} \pmb{n}_{\nu} = \pmb{n},$$
and from (\ref{alpha}) and (\ref{theta}) we deduce that 
\begin{equation}\label{5}
0 < |\pmb{\theta} - \pmb{\theta}_{\nu}| \ll \dfrac{a_{\nu}}{q_{\nu}}
\end{equation}
and 
\begin{equation}\label{angle}
{\rm angle} (\pmb{n},\pmb{n}_{\nu}) \ll \dfrac{b_{\nu}}{a_{\nu}}
\end{equation}
It is clear that $\pmb{\theta} \not\in \Q^2.$
A slight modification\footnote{A similar procedure was explained in \cite{NGM1}. There author provides the linear independence of coordinates of the limit vector by "going away from all rational subspaces" (the beginning of proof of Theorem 1 in the case $k=1$, p.132 and the beginning of section 5, p.146).} of the procedure of choosing vectors $\pmb{z}_{\nu}$ ensures the condition that $1,\theta_1,\theta_2$ are linear independent over $\Z$. Define $\pi = \{\pmb{x} \in \R^3: (\pmb{x},\pmb{n}) =0\}$. Then $\pmb{\theta} \in \pi$ by continuity and we can assume that $\pmb{n} \not\in \Q^3.$
\section{Winning property}
Consider  one-dimensional affine subspaces  $$\mathcal{P}_{\nu} = \{(x_1,x_2) \in \R^2 : (1,x_1,x_2) \in \pi_{\nu}\} \subset \R^2$$ and  $$\mathcal{P} = \{(x_1,x_2) \in \R^2 : (1,x_1,x_2) \in \pi\} \subset \R^2,$$ where $\pi$ was defined in the end of the previous section. 
Let
$$ B_1(\pmb{\theta}) = \{ \pmb{\xi} \in \R^2 : {\rm dist} (\pmb{\xi},\pmb{\theta}) < 1\}. $$
We will show that for any $\pmb{\eta} = (\eta_1, \eta_2) \in \mathcal{P} \cap B_1 (\pmb{\theta})$ there exists infinitely many solutions of the inequality
$$
\max_{i=1,2} 
 || \theta_{i} x- \eta_{i}|| < \dfrac{\omega(x)}{x}
$$
in integer $x.$
Denote by $\pmb{\eta}_{\nu} = (\eta_{1 \nu},\eta_{2 \nu})$ the orthogonal projection of $\pmb{\eta}$ onto $\mathcal{P}_{\nu}$.
From \eqref{angle} we see that
\begin{equation}\label{4}
|\pmb{\eta} -\pmb{\eta}_{\nu}| \ll \dfrac{b_{\nu}}{a_{\nu}}
\end{equation}
 
For any $\pmb{\eta}_{\nu} = (\eta_{1 \nu}, \eta_{2 \nu}) \in \PP_{\nu}$
 the planar domain $\overline{\pmb{\eta}_{\nu}}+ \Pi_{\nu}, \overline{\pmb{\eta}_{\nu}} = (1, -\eta_{1 \nu}, -\eta_{2 \nu})$ contains an integer point
$\pmb{y}_{\nu} = (x_{\nu},y_{1 \nu},y_{2 \nu}) \in \Lambda_{\nu}$.
It is clear that 
\begin{equation}\label{12}
|x_{\nu}|\ll q_{\nu}
\end{equation}
and 
\begin{equation}\label{13}
\max_{i=1,2} 
 | \theta_{i \nu} x_{\nu}- \eta_{i \nu}-y_{i \nu}|  \ll \dfrac{d_{\nu}}{q_{\nu}}.
\end{equation}
By \eqref{5},\eqref{4},\eqref{12},\eqref{13} we have
\begin{equation*}
\max_{i=1,2} 
 || \theta_{i} x_{\nu}- \eta_{i}|| \le
  |x_{\nu}| \max_{i=1,2}  |\theta_i - \theta_{i \nu}| +  \max_{i=1,2} ||\theta_{i \nu} x_{\nu}- \eta_{i \nu}|| + \max_{i=1,2} |\eta_i - \eta_{i \nu}| 
\end{equation*}  
 $$
 \ll a_{\nu} + \dfrac{d_{\nu}}{q_{\nu}} + \dfrac{b_{\nu}}{a_{\nu}} \ll \dfrac{d_{\nu}}{q_{\nu}}
$$
in the last inequality we use \eqref{b}.
By \eqref{7} we have 
$$
\max_{i=1,2} 
 || \theta_{i} x_{\nu}- \eta_{i}|| \le \dfrac{\omega(q_{\nu})}{q_{\nu}}
$$
for large $\nu.$
As $\pmb{\overline{\eta}} \in \pi$ and $\pmb{y}_{\nu} \in \pi_{\nu}$, 
$\max_{i=1,2} 
 || \theta_{i} x_{\nu}- \eta_{i}|| \neq 0$ infinitely often (in fact for all large $\nu$.) 


\begin{thebibliography}{99}
\bibitem{An} 
J. An, \emph{Two-dimensional badly approximable vectors and Schmidt's game},  Duke Math. J. 165, no. 2 (2016), 267-284

\bibitem{BMS}
P. Bengoechea, N. Moshchevitin, N. Stepanova \emph{A note on badly approximable linear forms on manifolds},  Mathematika, 63(2), 587-601. 

\bibitem{Cas}
J.W.S. Cassels, \,\, An introduction to Diophantine approximations, Cambridge Univ. Press, 1957.


 \bibitem{Dremov}

    V. A. Dremov, \emph{On domains of $(\alpha, \beta)$-winnability}, Dokl. Akad. Nauk 384 (2002), no. 3, 304-307 (Russian). 


\bibitem{ET} 
M. Einsiedler, J. Tseng, \emph{Badly approximable systems of affine forms, fractals, and Schmidt games}, J. Reine Angew. Math. 660 (2011), 83-97.



\bibitem{FSU}
L. Fishman, D. Simmons, M. Urba\'{n}ski, \emph{Diophantine approximation and the geometry of limit sets in
Gromov hyperbolic metric spaces}, http://arxiv.org/abs/1301.5630, preprint 2013, to appear in Mem. Amer. Math.
Soc.

\bibitem{HM} S. Harrap, N. Moshchevitin, \emph{A note on weighted badly approximable linear forms}, to appear in Glasgow Mathematical Journal.


\bibitem{J1}
V. Jarn\'{i}k,\emph{ O linea\'{r}n\'{i}ch nehomogenn\'{i}ch
  diofantick\'{y}ch aproximac\'{i}ch (on linear inhomogeneous Diophantine approximations)},
   Rozpravy II. T\v{r}\'{i}dy \v{C}esk\'{e} 
   Akad. 51 (1941), no. 29, 21. MR 0021015  
 
\bibitem{J2} 
 V. Jarn\'{i}k,\emph{ Sur les approximations diophantiques lin\'{e}aires non homog\'{e}nes}, Acad. Tch\'{e}que Sci. Bull. Int. Cl. Sci. Math. Nat. 47 (1946), 145-160 (1950).
  
 
 
\bibitem{KH1}
A. Khintchine, \emph{\"Uber eine Klasse linearer diophantischer Approximationen}, Rend. Circ. Mat. Palermo 50 (1926), 170-195.

\bibitem{KH2}
A. Khintchine, 
\emph{\"Uber die angen\"{a}herte Aufl\"{o}sung linearer Gleichungen in ganzen Zahlen},
 Acta Arith. 2 (1937), 161-172. 
\bibitem{KH3}   
    A.  Khintchine, \emph{Regular systems of linear equations and a general problem of Chebyshev}, Izvestiya Akad. Nauk SSSR. Ser. Mat. 12 (1948), 249-258 (Russian).   
 
 
 \bibitem{CH}
A. Khintchine, \emph{Continued fractions}, With a preface by B. V. Gnedenko. Translated from the third
(1961) Russian edition. Reprint of the 1964 translation. Dover Publications, Inc., Mineola, NY, 1997
 
\bibitem{WK}
D. Kleinbock, B. Weiss, \emph{Modified Schmidt games and Diophantine approximation with weights}, Advances in Math. 223 (2010), 1276-1298.


\bibitem{Mcm}
C. McMullen, \emph{Winning sets, quasiconformal maps and Diophantine approximation}, to appear in Geom. Funct. Anal. 20 (2010), 726-740.



\bibitem{NGM}
N. Moshchevitin, \emph{A note on badly approximable affine forms and winning sets}, Mosc. Math. J. 11 (2011), no. 1, 129-137.



\bibitem{NGM1}
N. Moshchevitin,
\emph{
    Proof of W.M. Schmidt's conjecture concerning successive minima of a lattice}
Journal of the London Mathematical Society, Volume 86, Issue 1 (2012), Pages 129-151

\bibitem{Sch1}
W. M. Schmidt, \emph{Badly approximable systems of linear forms}, J. Number Theory 1 (1969), 139-154. 

\bibitem{Sch2}
 W. M. Schmidt, \emph{Diophantine approximation}, Lecture Notes in Math., vol. 785, Springer-Verlag, Berlin 1980.
 
\bibitem{Sch3}
  W. M. Schmidt, \emph{On badly approximable numbers and certain games}, Trans. Amer. Math. Soc. 123 (1966), 178-199.
\bibitem{T} J. Tseng, \emph{Badly approximable affine forms and Schmidt games}, J. Number Theory 129 (2009), 3020-3025.

  
 
 
 
\end{thebibliography}
\end{document}